\documentclass{ifacconf}

\usepackage{natbib}            % you should have natbib.sty
\usepackage{graphicx}          % Include this line if your
\usepackage[T1]{fontenc}
\usepackage{subfig}
\usepackage{amssymb}
\usepackage{amsmath}
\usepackage{dcolumn}
\usepackage{bm}
\usepackage{color}
\usepackage{algorithm}
\usepackage{algorithmicx}
\usepackage{algpseudocode}
\usepackage{nomencl}
\usepackage{multirow}
\usepackage{slashbox}

\newcommand{\R}{\mathbb{R}}
\newcommand{\E}{\mathbb{E}}
\newcommand{\noise}{\sigma^2}
\newcommand{\hyper}{\gamma}
\newcommand{\bxi}{\boldsymbol{\xi}}

\newcommand{\bb}{\mathbf{b}}
\newcommand{\bc}{\mathbf{c}}
\newcommand{\be}{\mathbf{e}}
\newcommand{\bff}{\mathbf{f}}
\newcommand{\bF}{\mathbf{F}}

\newcommand{\br}{\boldsymbol{\gamma}}
\newcommand{\bx}{\mathbf{x}}
\newcommand{\bu}{\mathbf{u}}
\newcommand{\bv}{\mathbf{v}}
\newcommand{\bw}{\mathbf{w}}
\newcommand{\by}{\mathbf{y}}
\newcommand{\bz}{\mathbf{z}}

\newcommand{\bI}{\mathbf{I}}
\newcommand{\bL}{\mathcal{L}}
\newcommand{\bN}{\mathcal{N}}
\newcommand{\bP}{\mathcal{P}}
\newcommand{\bU}{\mathbf{\Theta}}
\newcommand{\bW}{\mathbf{W}}
\newcommand{\btheta}{\boldsymbol{\theta}}

\newcommand{\dic}{\mathbf{A}}
\newcommand{\dicc}{\mathbf{A}}

\newcommand{\bGamma}{\mathbf{\Gamma}}
\newcommand{\bXi}{\boldsymbol{\xi}}
\newcommand{\pan}{\noise\mathbf{I+{\dic}\bGamma}{{\dic}}^{\mathrm{T}}}
\newcommand{\pank}{\noise\mathbf{I+{\dic}\bGamma^k}{{\dic}}^{\mathrm{T}}}

\newcommand{\prior}{\varphi}
\newcommand{\define}{\triangleq}

\newcommand{\diag}{\operatorname{diag}}

\DeclareMathOperator*{\argmin}{argmin}
\newtheorem{assumption}{Assumption}
\newtheorem{remark}{Remark}
\begin{document}

\begin{frontmatter}

\title{Distributed Reconstruction of Nonlinear Networks: An ADMM Approach\thanksref{footnoteinfo}} % Title, preferably not more than 10 words.

\thanks[footnoteinfo]{W. Pan  gratefully acknowledge the support of Microsoft Research through
the PhD Scholarship Program. A. Sootla and G.-B. Stan acknowledge the support of EPSRC through the project EP/J014214/1 and the EPSRC Science and Innovation Award EP/G036004/1.}

\author{Wei Pan,}
\author{Aivar Sootla,}
\author{and Guy-Bart Stan}

\address{Centre for Synthetic Biology and
Innovation and the Department of Bioengineering, Imperial College London, United Kingdom; e-mail:  \{w.pan11, a.sootla, g.stan\}@ imperial.ac.uk}

\begin{abstract}                          % Abstract of not more than 250 words.
In this paper, we present a distributed algorithm for the reconstruction of large-scale nonlinear networks. In particular, we focus on the identification from time-series data of the nonlinear functional forms and associated parameters of large-scale nonlinear networks. In (\cite{pantac}), a nonlinear network reconstruction problem was formulated as a nonconvex optimisation problem based on the combination of a marginal likelihood maximisation procedure with sparsity inducing priors. Using a convex-concave procedure (CCCP), an iterative reweighted lasso algorithm was derived to solve the initial nonconvex optimisation problem. By exploiting the structure of the objective function of this reweighted lasso algorithm, a distributed algorithm can be designed. To this end, we apply the alternating direction method of multipliers (ADMM) to decompose the original problem into several subproblems. To illustrate the effectiveness of the proposed methods, we use our approach to identify a network of interconnected Kuramoto oscillators with different network sizes (500$\sim$100,000 nodes).
\end{abstract}

\end{frontmatter}

\section{Introduction}

%Reconstruction of large-scale networks for systems and processes described by discrete-time nonlinear models (also known as finite difference nonlinear models) from time-series data is very important in many different fields such as systems/synthetic biology, econometrics, finance, chemical engineering, social networks, etc. Yet techniques for reconstruction still remain challenging for system identification, (\cite{ljung2011four}). Graphical modelling approaches are widely used to learn dynamic models. Causal relations among nodes in networks can be encoded in a graphical way via ``Dynamic Bayesian Networks'', (\cite{kollar2009probabilistic,pearl1988probabilistic,spirtes2000causation}).
%``Causality'' of the dynamics of one node in a network comes from two sources: one is its own past history, the other one is the causal influence of the other state variables in the network and their associated past history. In the context of graphical modelling approaches, the dynamics are typically assumed to be linear. Under a linear dynamics assumption, graphical modelling approaches can be used to construct linear autoregressive (AR) models from limited time-series data, (\cite{bach2004learning,barber2010graphical}).
%However, such approaches typically fail to represent nonlinear dynamical systems.

The importance of reconstructing nonlinear systems and its associated difficulties are widely recognised (\cite{ljung2011four}).
Reconstruction methods focus on specific system classes such as those described by Wiener and Volterra series or nonlinear auto-regressive with exogenous inputs (NARX) models to name just a few examples (see \cite{ljung1999system} and references therein). However, nonlinear systems can be described by other functional forms. One of the most important and challenging problems in nonlinear network reconstruction is nonlinear structure identification (\cite{sjoberg1995nonlinear}). Nonlinear functional forms can be typically expanded as sums of terms belonging to a family of parameterised functions (see Sec.~5.4, \cite{ljung1999system}). A usual approach to identifying a nonlinear black-box model is to search amongst a set of possible nonlinear terms (e.g., basis functions) for a parsimonious description coherent with the available data set \cite{haber1990structure}).

In this paper, nonlinear systems are represented in a general state-space form. The framework
we develop uses some a priori knowledge of the type of system we want to reconstruct, i.e., we
consider a set of candidate dictionary functions appropriate for the specific type of systems from which the data have been collected (e.g., biological, chemical, mechanical, or electrical system).
We assume the measured time-series data and a set of candidate dictionary functions are given.
Our main objective is to identify the most parsimonious representation that explains the collected time-series data at best. Generally, we cast the reconstruction problem into a sparse signal recovery problem (\cite{Candes2005decoding,donoho2006compressed}).
In \cite{pantac}, the nonlinear network reconstruction problem was casted as a nonconvex optimisation problem which was shown to be efficiently solvable using a centralised reweighted lasso algorithm. Nonlinear reconstruction problems solved by centralised reconstruction methods have typically a relatively small size.

Social networks, communication networks and biological networks are typically very large (e.g. more than 100,000 nodes) and the data set collected from them is therefore quite ``big''. Typically, centralised reconstruction algorithms cannot handle such problems due to their associated very large memory and computational requirements.
Here, we will apply the alternating direction method of multipliers (ADMM) to split the centralised problem into several subproblems with each subproblem solving a weighed lasso problem independently. This approach has the advantage that memory and computational requirements can be both reduced in comparison to generic centralised solvers.% such as CVX (\cite{grant2008cvx}) or YALMIP (\cite{lofberg2004yalmip}).

ADMM is a powerful algorithm for solving structured convex optimization problems. The ADMM method was introduced for optimisation in the 1970's and is closely related to many other optimisation algorithms including Bregman iterative algorithms, Douglas-Rachford splitting, and proximal point methods (see (\cite{boyd2011distributed}) and references therein). ADMM has been shown to have strong convergence properties and to be useful for solving by decomposition large optimisation problems which cannot be handled by generic optimization solvers. ADMM has been applied in many areas, such as filtering (\cite{wahlberg2012admm}), image processing (\cite{figueiredo2010restoration}) as well as large-scale problems in statistics and machine learning (\cite{boyd2011distributed}).

The paper is organised as follows.
In Section \ref{sec:model}, we formulate the nonlinear network reconstruction problem we consider in this paper.
In Section~\ref{sec:bayesian}, we re-interpret this reconstruction problem from a Bayesian point of view.
In Section~\ref{sec:cccp}, we derive an iterative reweighted $\ell_1$ lasso algorithm to solve the nonconvex optimisation problem based on concave-convex procedure.
In Section~\ref{sec:admm}, we review ADMM, apply it to our optimisation problem, and derive a distributed algorithm.
In Section~\ref{sec:example}, we apply our method to the reconstruction of networks of interconnected Kuramoto oscillators.
Finally, in Section \ref{sec:conclusion}, we conclude and discuss several future problems.

%%%%%%%%%%%%%%%%%%%%%%%%%%%%%%%%%%%%%%%%%%%%%%%%%%%%%%%%%%%%%%%%%%%%%%%%%%%%%%%%%%%%
\section{PROBLEM FORMULATION}\label{sec:model}

\subsection{Nonlinear Dynamical Systems}
We consider dynamical systems described by multi-input multi-output (MIMO) nonlinear discrete-time equations with additive noise:
\begin{eqnarray}
\bx(t_{k+1})=\bF(\bx(t_{k}),\bu(t_{k}))+\bxi(t_{k}),
\label{model:vector}
\end{eqnarray}
where $\bx=[x_{1},\ldots,x_{n_\bx}]^{\mathrm{T}}\in {\R}^{n_\bx}$ denotes the state vector;
$\bu=[u_{1},\ldots,u_{n_\bu}]^{\mathrm{T}}\in {\R}^{n_\bu}$ denotes the input vector;
$\bF(\cdot)\define \left[\bF_1(\cdot),\ldots,\bF_{n_\bx}(\cdot)\right]^T: \mathbb{R}^{n_\bx+n_\bu}\rightarrow \mathbb{R}^{n_\bx}$,
and $\bxi= [\xi_{1},\xi_{2},\ldots,\xi_{n_\bx}]^{\mathrm{T}} \in{\R}^{n_\bx}$ is assumed to
be a zero-mean Gaussian white noise vector with constant positive covariance matrix $\mathbf{\Xi}$, i.e., $\bxi(t_{k})\thicksim\bN(\mathbf{0}, \mathbf{\Xi})$.
Since this description covers most of the discrete-time nonlinear dynamical systems with infinite number of possible functional forms for $\bF(\cdot)$, we confine the scope and assume system (\ref{model:vector}) satisfies the following assumptions:
\begin{assumption}
\label{assumption:observable}
The system (\ref{model:vector}) is fully measurable, i.e., all the state variables $x_i$ can be measured and there are no hidden variables.
\end{assumption}
\begin{assumption}
\label{assumption:dictionary}
The function terms $\bF(\bx(t_{k}),\bu(t_{k}))$ in (\ref{model:vector}) are smooth and can be represented as a linear combinations of several dictionary functions, (see Sec.~5.4 in \cite{ljung1999system}).
\end{assumption}

%\emph{Problem}. Based on the assumptions above, the problem that we are interested in solving is the
%following: given the time-series data of $\bx$, what's the form of the nonlinear function $\bF$?

\subsection{Construction of Dictionary Functions}

Depending on the field for which the dynamical model needs to be built, only a few typical nonlinearities specific to this field need to be considered. For example, the class of models that arise from genetic regulatory networks (GRN) typically involves nonlinearities that capture fundamental biochemical kinetic laws, e.g., first-order degradation functions, mass-action kinetics, Hill and Michaelis-Menten functions, which are confined to either polynomial or rational functions. In what follows we gather the set of \emph{all} candidate/possible  dictionary functions that we want to consider for reconstruction.
Consider state variable $x_{i}$, $i=1,\ldots,n_{\bx}$. Under Assumption \ref{assumption:dictionary}, the function terms for state $i$ can be written as:
\begin{equation}
\begin{aligned}
\bF_{i}(\bx(t_{k}),\bu(t_{k}))&=\sum_{s=1}^{N_i}w_{is}\bff_{is}(\bx(t_{k}),\bu(t_{k})),  \\
&=\bw_i^{\mathrm{T}}\bff_i(\bx(t_{k}),\bu(t_{k})),
\label{eq:expansion}
\end{aligned}
\end{equation}
where $\bw \in {\mathbb{R}^{N_i}}$ and $\bff_{i}: \mathbb{R}^{n_\bx+n_\bu}\rightarrow \mathbb{R}^{N_i}$ are dictionary functions that are assumed to govern the dynamics.  $\bff_{i}(\bx(t_{k}),\bu(t_{k}))$ can be monomial, polynomial, constant or any other functional form such as rational, exponential, trigonometric etc.

Taking the transpose of both sides of (\ref{eq:expansion}), we obtain
\begin{equation}
\begin{aligned}
x_i(t_{k+1})=\bff_i^T(\bx(t_{k}),\bu(t_{k}))\bw_i+\xi_i(t_{k}), \ i=1,\ldots,n_{\bx},
\label{model:transpose}
\end{aligned}
\end{equation}
where $\xi_i(t_{k})$ is assumed to be i.i.d. Gaussian distributed: $\xi_i(t_{k})\thicksim\bN(0, \noise_i)$, with
$
\E(\xi_i (t_{p}))=0, \ \E(\xi_i(t_{p})\xi_i(t_{q}))=\noise_i\delta_{pq},
\label{model:noise}
$
with
$
\delta _{pq}=\left\{
\begin{array}{ll}
1, & p=q, \\
0, & p\neq q
\end{array}
\right.
$.
If $M+1$ data samples including the initial value satisfying (\ref{model:transpose}) can be obtained from the system of interest, the system in (\ref{model:transpose}) can be written as
\begin{equation}
\begin{aligned}
\mathbf{y}_i=\dicc_i\bw_i+\bXi_i, \ i=1,\ldots,n_{\bx}.
\label{problem:expand}
\end{aligned}
\end{equation}
with
{\small
\begin{equation*}
\begin{aligned}
\by_i &\define
\left[x_i(t_{1}),\ldots,x_i(t_{M})\right]^{\mathrm{T}}\in {\mathbb{R}}^{M},\notag\\
\mathbf{\dicc}_i &\define
\left[
\begin{array}{ccc}
f_{i1}(\bx(t_{0}),\bu(t_{0})) & \ldots  & f_{iN_i}(\bx(t_{0}),\bu(t_{0})) \\
\vdots  & \ddots  & \vdots  \\
f_{i1}(\bx(t_{M-1}),\bu(t_{M-1})) & \ldots  & f_{iN_i}(\bx(t_{M-1}),\bu(t_{M-1}))
\end{array}
\right] \\
&\in {\mathbb{R}}^{M\times N_i}, \notag \\
\bw_i &\define \left[w_{i1},\ldots,w_{iN_i}\right]^T \in {\mathbb{R}}^{N_i},\notag \\
\bXi_i &\define
\left[\xi_i(t_{0}),\ldots,\xi_i(t_{M-1})\right]^{\mathrm{T}}\in {\mathbb{R}}^{M}.
\label{formulation:0}
\end{aligned}
\end{equation*}
}
\hspace{-.05 in}Since the $n_{\bx}$ linear regression problems in (\ref{problem:expand}) are independent, for simplicity of notation, we omit the subscript $i$ in (\ref{problem:expand}) and write
\begin{equation}
\begin{aligned}
\mathbf{y}=\dic \bw+\bXi.
\label{problem}
\end{aligned}
\end{equation}
The problem is thus to find $\bw$ given the measured noisy data stored in $\mathbf{y}$.

%This a typical linear regression problem that can be solved using standard least square approaches, provided that the structure of the nonlinearities in the model are known, i.e., provided that $\dicc_i$ is known.

\subsection{Discussion on relaxations of solutions}
\label{sec:compressive sensing}

Considering the network reconstruction problem in practice, several problems arise with respect to  (\ref{problem}).
Firstly, a low number of time-series measurements will render the linear regression in (\ref{problem}) under-determined.
Secondly, the number of columns of the dictionary matrix might be very large, due to the potential introduction of non-relevant and/or non-independent dictionary functions in $\dic$.
 As a result, the sparsest solution will be favoured due to the model selection criterions such as Akaike information criterion, (\cite{akaike1974new}) or Bayesian information criterion, (\cite{schwarz1978estimating}).
%
%\guy{MIGHT NEED TO BE FIXED (I HAVE TRIED BY DELETING WHAT I FELT WAS INCOMPREHENSIBLE)}
%% As a result, the sparsest solution will be favoured due to the model selection criterion such as
%
%% the true solution $\bw$ to \eqref{problem:expand} is typically going to be sparse.
However, finding the sparsest solution is NP-hard. Classically, a lasso algorithm is typically used as a relaxation (\cite{tibshirani1996regression}) to this NP-hard optimisation problem.
Lasso usually works well when the dictionary matrix has certain properties such as the \textit{restricted
isometry property} (RIP), (\cite{Candes2005decoding} or the incoherence property, (\cite{donoho2003optimally}).
These properties basically state that two or more of the columns of dictionary matrix cannot be co-linear or close to be co-linear.
Unfortunately, such properties are hardly guaranteed in typical network reconstruction problems.
%
%\guy{(I HAVE TRIED BY DELETING WHAT I FELT WAS INCOMPREHENSIBLE)}
%% $\ell_{1}$ relaxation based algorithms
%% do not yield the optimal solution to the network reconstruction problem.
In the following, we shall introduce, from a probabilistic viewpoint, how a Bayesian treatment can alleviate these RIP or incoherence requirements (\cite{tipping2001sparse,seeger2010variational}).

\section{BAYESIAN VIEWPOINT}
\label{sec:bayesian}
\subsection{Sparsity Inducing Prior}
Bayesian modelling treats all unknowns as stochastic variables with certain probability distributions, (\cite{bishop2006pattern}). For $\by=\dic \bw+\bXi$, it is assumed that the stochastic variables in $\bXi$ are i.i.d. Gaussian distributed with $\bXi\thicksim\bN(\mathbf{0}, \noise\bI)$.
%We further define the precision or inverse-variance as $\beta=1/\sigma^{2}$.
In such case, the likelihood of the data given $\bw$ is
\begin{equation}
\begin{aligned}
\bP(\by|{\bw})
={\mathcal{N}}(\by|{\dic} {\bw},\sigma^{2}\bI)
\propto\exp \left[ -\frac{1}{2\sigma^{2}} \| \dic\bw-\by\|_2^{2}\right].
\label{likelihood}
\end{aligned}
\end{equation}
Given the likelihood function in (\ref{likelihood}) and specifying a prior
$\bP(\bw)=\prod_{j}\bP(w_j)$.
%we compute the \emph{posterior distribution} over $\bw$ via Bayes' rule:
%$\bP(\bw|{\by})\propto {\bP(\by|\bw)\bP(\bw)}$.
We further define a prior distribution $\bP(\bw)$ as
$
\bP(\bw)\propto \exp \left[-\frac{1}{2}g(\bw)\right]=\exp \left[-\frac{1}{2}\sum_{j}g(w_j)\right],
$
where $g(w_j)$ is a given function of $w_j$.
To enforce sparsity on $\bw$, the function $g(\bw)$ is usually chosen as a concave, non-decreasing function of $|\bw|$.
Such penalty functions can be realised by adding sparsity inducing priors on $\bw$. Such sparsity inducing priors include  Laplace $\propto \exp (-\gamma\sum_{j}|{w}_j|)$, where $\gamma>0$, Student's \emph{t}$\propto(b+w_j^2/2)^{-(a+\frac{1}{2})}$, where $a, b>0$, and others, (\cite{palmer2005variational}).

Typically,
$\bP(\bw|\by)$ is approximated by a Gaussian distribution from which efficient algorithms exist  (\cite{bishop2006pattern}). To this end, we may consider \emph{super-Gaussian} priors, which yield a lower bound for the priors $\bP(w_j)$ (\cite{palmer2005variational}). More specifically, if we define $\boldsymbol{\hyper} \define \left[\hyper_1, \ldots, \hyper_N\right]^T \in \R^N_{+}$, we can represent the prior in the following relaxed (variational) form:
\begin{equation}
\begin{aligned}
\bP(\bw)=\prod_{j=1}^{n}\bP(w_j),\ \bP(w_j)=\max_{\hyper_j>0}\bN(w_j|0,\hyper_j)\prior(\hyper_j),
\label{prior2}
\end{aligned}
\end{equation}
where $\prior(\hyper_j)$ is a nonnegative function which is treated as a hyperprior with $\hyper_j$ being its associated hyperparameters. Throughout, we call $\prior(\hyper_j)$ the ``\emph{potential function}''. This Gaussian relaxation is possible if and only if $\log \bP(\sqrt{w_j})$ is
concave on $(0,\infty)$.
Hereafter, we adopt Student's \emph{t} prior and the corresponding $\prior(\hyper_j)=1$.

\subsection{Marginal Likelihood Maximisation}
For a fixed $\br=\left[{\hyper}_{1},\ldots, {\hyper}_{N}\right]$,
we define a relaxed prior which is a joint probability distribution over $\bw$ and ${\br}$
\begin{equation}
\begin{aligned}
\bP(\bw; \br)
=\prod_{j}\bN(w_j|0,{\hyper}_j)\prior({\hyper}_j)
=\bP(\bw|\br)\bP(\br) \le \bP(\bw),
\label{prior-relaxation}
\end{aligned}
\end{equation}
where $\bP(\bw|\br)\define \prod_{j}\bN(w_j|0,{\hyper}_j),
\bP(\br) \define \prod_{j} \prior({\hyper}_j).$
Since is  $\bP(\by|\bw)$ is Gaussian in (\ref{likelihood}), we can get a relaxed posterior which is also Gaussian.

Now the key question is how to choose the most appropriate $\br=\hat{\br}=\left[\hat{\hyper}_{1},\ldots, \hat{\hyper}_{N}\right]$ to maximise $\prod_{j}\bN(w_j|0,{\hyper}_j)\prior({\hyper}_j)$ such that $\bP(\bw|\by,\hat{\br})$ can be a ``good'' relaxation to $\bP(\bw|\by)$.
Using the product rule for probabilities, we can write the full posterior
$
\bP(\bw, \br|\by)
\propto \bN(\mathbf{m}_{\bw},\mathbf{\Sigma}_{\bw}) \times \frac{\bP(\by|\br)\bP(\br)}{\bP(\by)}.
$
Since $\bP(\by)$ is independent of $\br$, the quantity
$$\bP(\by|\br)\bP(\br)=\int \bP(\by|\bw)\bP(\bw|\br)\bP(\br)d\bw$$ is the prime target for variational methods (\cite{wainwright2008graphical}). This quantity is  known as evidence or marginal likelihood.
A good way of selecting $\hat{\br}$ is to choose it as the minimiser of the sum of the misaligned probability mass, e.g.,
\begin{equation}
\begin{aligned}
\hat{\br}&=  \argmin\limits_{\br\geq\mathbf{0}} \int \bP(\by|\bw)\left|\bP(\bw)-\bP(\bw;\br)\right|d\bw \notag \\
&= \arg \max\limits_{\br\geq\mathbf{0}} \int \bP(\by|\bw)\prod_{j=1}^{n}\bN(w_j|0,\hyper_j)\prior(\hyper_j)d\bw.
\label{mass}
\end{aligned}
\end{equation}
The second equality is  a consequence of $\bP(\bw;\br)\le\bP(\bw)$ (see (\ref{prior-relaxation})).
The procedure in (\ref{mass}) is referred to as  evidence/marginal likelihood maximisation or type-II maximum likelihood, (\cite{tipping2001sparse,seeger2010variational, wipf2011latent, seeger2011large}).
It means that the marginal likelihood can be maximised by selecting the most probable hyperparameters able to explain the observed data.
Once $\hat{\br}$ is computed, an estimate of the unknown weights can be obtained by setting $\hat{\bw}$ to the posterior mean %(\ref{mean}):
\begin{equation}
\begin{aligned}
\hat{\bw}=\E(\bw|\by;\hat{\br})=\hat{\bGamma} \dic^T (\noise\mathbf{I+{\dic}\hat{\bGamma}}{{\dic}}^{\mathrm{T}})^{-1} \by.
\label{wsbl}
\end{aligned}
\end{equation}
with $\hat{\bGamma}\define \diag[\hat{\br}]$.
If an algorithm can be proposed to compute $\hat{\br}$ in (\ref{mass}), we can obtain an estimation of the posterior mean $\hat{\bw}$ in (\ref{wsbl}).

In (\cite{tipping2001sparse,wipf2011latent}), it is shown that $\hat{\bGamma}$ can be obtained by minimising the following nonconvex cost function
\begin{equation}
\begin{aligned}
{\bL_{\br}}\left( \br \right) =\log \left\vert \pan \right\vert +\by^{\mathrm{T}}(\pan)^{-1}\by.
\label{hypercostfunction}
\end{aligned}
\end{equation}
However, such formulation does not allow to incorporate convex constraints on $\bw$, which are typically very useful in network reconstruction.
In (\cite{wipf2011latent,pantac}), it is shown that the dual cost function in $\bw$-space has the following form
\begin{equation}
\label{w-nonconvex}
\min_{\br\geq \mathbf{0}, \bw} \|\dic\bw-\by\|_{2}^{2}+\noise \bw^T\bGamma^{-1}\bw+\log|\pan|,
\end{equation}
To ease notation and avoid interrupting the flow, we will derive the algorithm without considering convex constraints in the sequel.

%%%%%%%%%%%%%%%%%%%%%%%%%%%%%%%%%%%%%%%%%%%%%%%%%%%%%%%%%%%%%%%%%%%%%%%%%%%%%%%%%%%%
\section{CONCAVE-CONVEX PROCEDURE}\label{sec:cccp}

The cost function in (\ref{w-nonconvex}) is convex in $\bw$ but nonconvex in $\bGamma$. In the next section, we show how this nonconvex optimisation problem can be formulated as a concave-convex procedure (CCCP) and then, show that finding the solution to CCCP is equivalent to solving iterative reweighted lasso problem. CCCP is another interpretation of the derivation of iterative reweighted lasso algorithm (\cite{wipf2010iterative}).

Let
\begin{equation}
\begin{aligned}
u(\bw, \br)&\define  \|\dic\bw-\by\|_{2}^{2}+\noise \sum_j \frac{w_j^2}{\hyper_j}, \notag \\
v(\br) &\define -\log|\pan| \notag
\end{aligned}
\end{equation}
Note that $u(\bw,\br)$ is jointly convex in $\bw$ and $\br$, $v(\br)$ is convex in $\br$.
Then the minimisation of the cost function (\ref{w-nonconvex}) can be formulated as a  concave-convex procedure
\begin{equation}
\begin{aligned}
\label{cccp}
\min_{\br\geq \mathbf{0}, \bw} u(\bw, \br)-v(\br)
\end{aligned}
\end{equation}
Since $v(\br) $ is differentiable over $\br$, the problem in (\ref{cccp}) can be transformed into the following iterative convex optimisation problem
\begin{equation}
\begin{aligned}
\left[\bw^{k+1},\br^{k+1}\right]
=\argmin\limits_{\br\geq \mathbf{0}, \bw} u(\bw, \br)-\nabla_{\br} v(\br^k)^T\br.
\label{cccp-3}
\end{aligned}
\end{equation}

Using basic principles in convex analysis, we then obtain the following analytic form for the negative gradient of $v(\br)$ at $\br$is:
\begin{equation}
\begin{aligned}
\boldsymbol{\alpha}^{k}
&=-\nabla_{\br} v(\br^k)^T \notag \\
&=-\nabla_{\br}\left( -\log |\pan| \right)|_{\br=\br^k}\notag \\
&= \diag \left[\dic^{\mathrm{T}}\left(\pank\right)^{-1}{\dic} \right],
\label{gammastarupdate}
\end{aligned}
\end{equation}

Then the iterative procedure (\ref{cccp-3}) can be formulated as
\begin{equation}
\begin{aligned}
&\left[\bw^{k+1},\br^{k+1}\right] =\\
&\argmin\limits_{\br\geq \mathbf{0}, \bw} \| \dic\bw-\by\|_{2}^{2}+\noise\sum_j \left(\frac{w_j^2}{\hyper_j}+\alpha_j^{k}\hyper_j\right).
\label{cccp-4}
\end{aligned}
\end{equation}
%%%%%%%%%%%%%%%%%%%%%%%%%%%%%%%%%%%%%%%%%%%%%%%%%%%%%%%%%%%%%%%%%%%%%%%%%%%%%%%%%%%%%

%We define the cost function to be minimised in (\ref{cccp-4}) as $\bL^k(\br,\bw)$.
The objective function in (\ref{cccp-4}) is jointly convex in $\bw$ and $\br$ and can be globally minimised by solving over $\br$ and then $\bw$. If $\bw$ is fixed, it gives
\begin{equation}
\begin{aligned}
\br^{k+1}=\argmin\limits_{\br\geq \mathbf{0}}\| \dic\bw-\by\|_{2}^{2}+\noise\sum_j \left(\frac{w_j^2}{\hyper_j}+\alpha_j^{k}\hyper_j\right). \label{1-gamma}
\end{aligned}
\phantom{\hspace{6cm}}
\end{equation}
We notice that in (\ref{1-gamma}), $\br^{k+1}$ has closed form solution $
\gamma_j^{k+1}=|w_j|/\sqrt{\alpha_{j}^k}$,
If $\gamma_j^{k+1}={|{w}_j|}/{\sqrt{\alpha_{j}^k}}$ is substituted into (\ref{1-gamma}), we get
\begin{equation}
\begin{aligned}
\bw^{k+1}
&= \argmin\limits_{\bw}\| \dic\bw-\by\|_{2}^{2}+\noise\sum_j \left(\frac{w_j^2}{\hyper_j^{k+1}}+\alpha_j^{k}\hyper_j^{k+1}\right) \notag\\
&=\argmin\limits_{\bw}\{\| \dic\bw-\by\|_{2}^{2}+2\noise \sum_{j=1}^{N}\sqrt{\alpha_{j}^k}|w_j| \}.
\label{cccp-5}
\end{aligned}
\end{equation}
We can then set
$
\gamma_j^{k+1}=|w_j^{k+1}|/\sqrt{\alpha_{j}^k}
\label{gammaupdate}
$, $\forall j$.  Then we update $\boldsymbol{\alpha}^{k+1}$ by (\ref{gammastarupdate}).
However, some of the estimated weights will be several orders of magnitude lower than the average ``energy'', e.g., $w_j^2 \ll \|\bw\|_2^2$. Thus a threshold needs to be defined \emph{a priori} to prune ``small'' weights at each iteration.
%It's proved in (\cite{wipf2010iterative} that the reweighted lasso in continued iteration can never do worse.
%Thereafter at each iteration we can record the index of the weights which are estimated to be zero, certain dictionary functions spanning the corresponding columns of $\dic$ can be pruned out for the next iteration.

We can now explain how the update of the parameters can be performed based on the above. Set the iteration count $k$ to zero and $\theta_j^0=1, \ \forall j$. At this stage, the optimisation is a typical lasso.
Then at the $k^{th}$ iteration, we initialise $\theta_j^{(k)}=\sqrt{\alpha^k_j}$, $\forall j$.
The above described procedure is summarised in Algorithm \ref{alg:reweightedlasso}.

%To this end, we can carry out the centralised optimisation in Algorithm \ref{alg:reweightedlasso} by using CVX (\cite{grant2008cvx}) or YALMIP (\cite{lofberg2004yalmip}), MATLAB packages for specifying and solving convex optimization problems. CVX or YALMIP calls generic SDP solvers (SDPT3 (\cite{toh1999sdpt3} or SeDuMi (\cite{sturm1999using}) to solve the problem. While these solvers are reliable for wide classes of optimisation problems, they are not customised for particular problem families, such as ours. Therefore, when the size of the network $n_{\bx}$ is very large or the number of candidate functions is very large or both, as we will illustrate in the examples later, the computational complexity is super expensive and the convergence is extremely slow.
 In the next section, we will reformulate the centralised optimisation in Algorithm \ref{alg:reweightedlasso} into distributed optimisation by ADMM.

\begin{algorithm}[!]
\caption{Reweighted lasso on $\bw$}
\label{alg:reweightedlasso}
\begin{algorithmic}[1]
\State Initialise $\theta_j^0=1, \ \forall j$%, $\dic^0=\dic$%,the weight pruning threshold $\epsilon$.
     \For  {$k=0, \ldots, k_{\max}$}
            \State Solve the weighted lasso problem
            \begin{equation}
            \bw^{k+1}=\argmin_{\bw}  \frac{1}{2}\| \dic\bw-\by\|_{2}^{2}+\noise\|\bU^{k}\bw\|_1; \label{weightedlasso}
            \end{equation}

%            \State Record the index $\mathbb{I}^k$ of nonzero entries in $\bw^{k+1}$, prune out the corresponding columns of $\dic$ indexed by  $\mathbb{I}^k$ to form $\dic^{k+1}$;

            \State Set
            $\bU^{(k)}\define \diag\left[\btheta^{(k)}\right]^{-1}$, $\bW^{(k)}\define \diag\left[|\bw^{(k)}|\right]$ and $
            \theta_j^{k+1}=\left[(\dic_j)^{\mathrm{T}}\left(\noise\mathbf{I}+\dic (\bU^k)^{-1} \bW^{k+1} (\dic)^T\right)^{-1}{\dic_j}\right]^{\frac{1}{2}}$

        \If   {A stopping criterion is satisfied}
            \State Break;
        \EndIf
     \EndFor
\end{algorithmic}
\end{algorithm}

\begin{remark}
\label{rmk:rl2}
There is one interesting finding from (\ref{1-gamma}). We can actually minimise over $\bw$ first then minimise over $\br$. Rather than only initialise $\boldsymbol{\alpha}^{0}$, we also initialise $\br^{0}$.  The iterative procedure is then a reweighted $\ell_2$ algorithm (\cite{wipf2010iterative})
\begin{equation}
\begin{aligned}
\bw^{k+1}
= \argmin\limits_{\bw}\| \dic\bw-\by\|_{2}^{2}+\noise\sum_j \frac{w_j^2}{\hyper_j^k}. \label{2-w}
\end{aligned}
\end{equation}
Then with the estimate $\bw^{k+1}$
\begin{equation}
\begin{aligned}
\br^{k+1}
= \argmin\limits_{\br\geq \mathbf{0}}\sum_j \left(\frac{(w_j^{k+1})^2}{\hyper_j}+\alpha_j^{k}\hyper_j\right). \label{2-gamma}
\end{aligned}
\end{equation}
$\br^{k+1}$ has a closed form solution $\br^{k+1}=|w_j^{k+1}|/\sqrt{\alpha_j^{k}}$. Then we can update $\boldsymbol{\alpha}^{k+1}$ as in (\ref{gammastarupdate}) and iterate until a stopping criterion is satisfied.

The $\bw$ update in (\ref{2-w}) is actually a reweighted-$\ell_2$ algorithm. $\br$ are the weighting vectors like $\btheta$ in reweighted lasso in Algorithm \ref{alg:reweightedlasso}. We also tested this reweighted $\ell_2$ algorithm for the examples in the sequel. However the heuristic convergence rate and reconstruction accuracy is not as good as reweighted lasso algorithm. We are still studying on it.

\end{remark}

%%%%%%%%%%%%%%%%%%%%%%%%%%%%%%%%%%%%%%%%%%%%%%%%%%%%%%%%%%%%%
%\section{distributed}
\section{ALTERNATING DIRECTION METHOD OF MULTIPLIERS (ADMM)}
\label{sec:admm}
In this section we give an overview of ADMM. We follow closely the development of (\cite{boyd2011distributed}
\subsection{ADMM}
ADMM is a numerical algorithm for solving optimisation problems such as
\begin{equation}
\label{w-nonconvex}
\begin{split}
\min_{\bw}\,\,\,\,\,\, & \,f(\bw)+g(\bz),\\
\mathrm{subject}\,\,\mathrm{to}\,\,\,\,\,\,\, &
P\bw+Q\bz=\bc,
\end{split}
\end{equation}
with variable $\bw \in \R^n$ and $\bz \in \R^m$, where $P\in \R^{p\times n}$, $Q \in \R^{p \times m}$, and $\bc \in \R^p$. We will assume that $f$ and $g$ are convex.

As the method of multipliers, we form the augmented Lagrangian
\begin{equation}
\begin{aligned}
L_{\rho}=&f(\bw)+g(\bz)+\bv^T(P\bw+Q\bz-c)+\\
&{\rho}/{2}\|P\bw+Q\bz-c\|_2^2.
\end{aligned}
\end{equation}
Defining the residual $r=P\bw+Q\bz-c$ and $\bu=(1/\rho)\bv$ as the scaled dual variable, we can express ADMM as
\begin{equation}
\begin{aligned}
\bw^{k+1}&:= \argmin\limits_{\bw}\left(f(\bw)+\frac{\rho}{2}\|P\bw+Q\bz^k-\bc+\bu^k\|_2^2\right)  \notag\\
\bz^{k+1}&:=\argmin\limits_{\bz}\left(g(\bz)+\frac{\rho}{2}\|P\bw^{k+1}+Q\bz-\bc+\bu^k\|_2^2\right)  \notag\\
\bu^{k+1}&:= \bu^k+P\bw^{k+1}+Q\bz^{k+1}-\bc.
\end{aligned}
\end{equation}

\subsubsection{Stopping criterion}
The primal and dual residuals at iteration $k$ are given by
$$\be_{primal}^k=\bw^k-\bz^k, \be_{dual}^k=-\rho(\bz^k-\bz^{k-1}).$$
We terminate the algorithm when the primal and dual residuals satisfy a stopping criterion:
$$\|\be_{primal}^k\|_2\leq\epsilon_{primal}, \|\be_{dual}^k\|_2 \leq \epsilon_{dual}.$$
Here, the tolerances $\epsilon_{primal}>0$ and $\epsilon_{dual}>0$ and can be set via an absolute plus relative criterion
\begin{equation}
\begin{aligned}
\epsilon_{primal}&= \sqrt{n}\epsilon_{abs}+ \epsilon_{rel}\max(\|\bw^k\|_2,\|\bz^k\|_2),\notag\\
\epsilon_{dual}&=  \sqrt{n}\epsilon_{abs}+ \epsilon_{rel}\rho\|\bu^k\|,
\end{aligned}
\end{equation}
where $\epsilon_{abs}$ and $\epsilon_{rel}$ are absolute and relative tolerances.
More details can be found in (\cite{boyd2011distributed}).

%%%%%%%%%%%%%%%%%%%%%%%%%%%%%%%%%%%%%%%%%%%%%%%%%%%%%%%
\subsection{Splitting across candidate functions}
In our setting, the number of candidate functions will be very large.
Therefore we partition across the candidate functions. Each subsystem can deal with its split of candidate functions independently then update the shared variables. The following are direct consequences of the so-called sharing problem in (\cite{boyd2011distributed}).

We partition the parameter vector $\bw$ as $\bw=(\bw_1,\ldots,\bw_n)$, with $\bw_i\in \R^{N_i}$, where $\sum_{i=1}^n N_i=N$. Partition the dictionary matrix $\dic$ as $\dic=[\dic_1,\ldots, \dic_n]$, with $\dic_i \in \R^{M \times N_i}$. Thus $\dic \bw=\sum_{i=1}^n \dic_i\bw_i$, i.e., $\dic_i\bw_i$ can be thought of as a `partial' prediction of $\by$ using only the candidate functions referenced in $\bw_i$.

Then the reweighted lasso problem (\ref{weightedlasso})
$$\min_{\bw} \frac{1}{2}\| \dic\bw-\by\|_{2}^{2}+\lambda\|\bU\bw\|_1$$
becomes
$$\min_{\bw}\frac{1}{2}\|\sum_{i=1}^{N}\dic_i\bw_i-\by\|_{2}^{2}+\lambda\sum_{i=1}^{N}\|\bU_i\bw_i\|_1.$$

Following the approach used for the sharing problem, we express the problem as
\begin{equation}
\begin{split}
\min\,\,\,\,\,\, & \, \frac{1}{2}\|\sum_{i=1}^{N}\bz_i-\by\|_{2}^{2}+\lambda\sum_{i=1}^{N}\|\bU_i\bw_i\|_1, \notag\\
\mathrm{subject}\,\,\mathrm{to}\,\,\,\,\,\,\, &
\, \dic_i\bw_i-\bz_i=0,\, \, i=1,\ldots, N,
\end{split}
\end{equation}
with new variables $\bz \in \R^m$. The derivation and simplification of ADMM also follows that for the sharing problem. The scaled form of ADMM is
{\small
\begin{align}
\bw_i^{k+1}&:= \argmin\limits_{\bw_i}\frac{\rho}{2}\|\dic_i\bw_i-\bz^k+\bu^k\|_2^2
+\lambda\|\bU_i\bw_i\|_1  \notag\\
\bz^{k+1}&:=\argmin\limits_{\bz} \frac{1}{2}\|\sum_{i=1}^{N}\bz_i-\by\|_{2}^{2}+\sum_{i=1}^N\frac{\rho}{2}\|\dic\bw^{k+1}-\bz^k+\bu^k\|_2^2\notag \\
\bu^{k+1}&:=\bu^k+\dic\bw^{k+1}-\bz^{k+1}.
\end{align}
}
\hspace{-.04 in}As in the discussion for the sharing problem, we carry out the $\bz$-update by first solving for the average
\begin{align}
\overline{\bz}^{k+1}&:=\argmin\limits_{\overline{\bz}} \|N\overline{\bz}-\by\|_2^2+\sum_{i=1}^N\frac{N\rho}{2}\|\overline{\bz}-\overline{\dic\bw}^{k+1}-\overline{\bu}^k\|_2^2
\notag \\
\bz^{k+1}&:=\overline{\bz}^{k+1}+\dic_i\bw_i^{k+1}+\bu^k-\overline{\dic\bw}^{k+1}-\overline{\bu}^k,
\end{align}
where $\overline{\dic\bw}^{k+1}=(1/N)\sum_{i=1}^N\dic_i\bw_i^{k+1}$. Substituting the last expression into the update for $\bu_i$, we find that
$$\bu_i^{k+1}=\overline{\dic\bw}^{k+1}+\overline{\bu}^k-\overline{\bz}^{k+1},$$
which shows that, as in the sharing problem, all the dual variables are equal. Using a single dual variable $\bu^k \in \R^m$, eliminating $\bz_i$, and define
$$\bb=\dic_i\bw_i^k+\overline{\bz}^k-\overline{\dic\bw}^k-\bu^k$$
we arrive at the Algorithm \ref{alg:splitting}
\begin{algorithm}[!]
\caption{ADMM for splitting across candidate functions}
\label{alg:splitting}
\begin{algorithmic}[1]
     \For  {$k=0, \ldots, k_{\max}$}
            \State Solve $$\bw_i^{k+1}= \argmin\limits_{\bw_i}\left(\frac{\rho}{2}\|\dic_i\bw_i-\bb\|_2^2
+\lambda\|\bU_i\bw_i\|_1\right);$$
            \State Compute $\overline{\bz}^{k+1}=\frac{1}{N+\rho}\left(\by+\rho\overline{\dic\bw}^{k+1}+\rho\bu^k\right)$;
            \State Update
            $\bu^{k+1}=\bu^k+\overline{\dic\bw}^{k+1}-\overline{\bz}^{k+1}$;
        \If   {A stopping criterion is satisfied}
            \State Break;
        \EndIf
     \EndFor
\end{algorithmic}
\end{algorithm}

%\begin{eqnarray}
%  \notag\\
%\overline{\bz}^{k+1}&:=&\argmin\limits_{\overline{\bz}} \left(\|N\overline{\bz}-\by\|_2^2+\sum_{i=1}^N\frac{N\rho}{2}\|\overline{\bz}-\overline{\dic\bw}^{k+1}-\bu^k\|_2^2
%\right)\notag \\
%&\,=& \frac{1}{N+\rho}\left(\by+\rho\overline{\dic\bw}^{k+1}+\rho\bu^k\right)\notag \\
%\bu^{k+1}&:=&\bu^k+\overline{\dic\bw}^{k+1}-\overline{\bz}^{k+1}.
%\end{eqnarray}

Each $\bw_i$-update is a lasso problem with $n_i$ variables, which can be solved using any lasso method.

In the $\bw_i$-update, we have $\bw_i^{k+1}=\mathbf{0}$ (meaning that none of the features in the $i$-th block are used) if and only if
\begin{equation}
\|\dic_i^T\bb\|_2
\leq \frac{\lambda}{\rho}
\end{equation}

The first step involves solving $N$ parallel weighted lasso (weighted $\ell_1$-regularised least squares) problems in $n_i$ variables each.
Between the first and second steps, we collect and sum the partial predictors $\dic_i\bw_i^{k+1}$ to form $\overline{\dic\bw}^{k+1}$. The second step is a single minimisation in $m$ variables, a quadratically regularised loss minimisation problem

\subsection{ADMM for Weighted Lasso}
%In the $\bw_i$-update, we need solve $n$ parallel weighted lasso (weighted $\ell_1$-regularised least squares) problems. This problem, however, as far as we are aware, there is no analytical formula. We carry out these updates by using CVX (\cite{grant2008cvx}) or YALMIP (\cite{lofberg2004yalmip}), Matlab packages for specifying and solving convex optimization problems. CVX or YALMIP calls generic SDP solvers (SDPT3 (\cite{toh1999sdpt3} or SeDuMi (\cite{sturm1999using}) to solve the problem. While these solvers are reliable for wide classes of optimisation problems, they are not customised for particular problem families, such as ours.

The weighted lasso problem can be solved using ADMM as well.
% In ADMM form, the problem can be written as
\begin{equation}
\begin{split}
\min\,\,\,\,\,\, & \, \frac{1}{2}\|\dic_i\bw_i-\bb\|_2^2
+\frac{\lambda}{\rho}\|\bU_i\bw_i\|_1, \notag\\
\mathrm{subject}\,\,\mathrm{to}\,\,\,\,\,\,\, &
\, \bU_i\bw_i-\hat{\bz}_i=0,\, \, i=1,\ldots, N.
\end{split}
\end{equation}
Define $\hat{\lambda}=\lambda/\rho$, it yields the ADMM algorithm
\begin{algorithm}[!]
\caption{ADMM for weighted lasso}
\label{alg:weightedlasso}
\begin{algorithmic}[1]
     \For  {$k=0, \ldots, k_{\max}$}
            \State Update    $\bw_i^{k+1}=(\dic_i^T\dic_i+\hat{\rho}\bU_i^T\bU_i)^{-1}(\dic_i^T\bb+\hat{\rho}\bU_i^T(\hat{\bz}_i-\hat{\bu}_i))$;
            \State Update $\hat{\bz}_i^{k+1}=S_{\hat{\lambda}/\hat{\rho}}(\bU_i\bw_i^{k+1}+\hat{\bu}_i^k)$;
            \State Update
            $\hat{\bu}_i^{k+1}= \hat{\bu}_i^k+\bU_i\bw_i^{k+1}-\hat{\bz}_i^{k+1}$;
        \If   {A stopping criterion is satisfied or when $k$ reaches a predefined iteration number $k_{\max}$}
            \State Break;
        \EndIf
     \EndFor
\end{algorithmic}
\end{algorithm}

 where $\hat{\rho}$ is the penalty parameter and the soft thresholding operator $S_{\hat{\lambda}/\hat{\rho}}$ is defined as
 $$S_{\hat{\lambda}/\hat{\rho}}(x)=\max(0, x-\hat{\lambda}/\hat{\rho})-\max(0, -x-\hat{\lambda}/\hat{\rho}).$$

\subsection{Algorithm For Nonlinear Network Reconstruction}
Now, we summarise the procedure for nonlinear network reconstruction in Algorithm \ref{alg:summary}.

\begin{algorithm}[!]
\caption{ADMM on $\bw$}
\label{alg:summary}
\begin{algorithmic}[1]
\State Initialisation \begin{enumerate}
\item Collect the time-series data, specify the candidate functions and construct the dictionary matrix;
%                        \item ADMM parameters $\lambda$, $\rho$, $\hat{\lambda}$ and $\hat{\rho}$; $\bU^{0}$
%                        \item Specify the splitting information
\item Partition the dictionary matrix $\dic$ as $\dic=[\dic_1,\ldots, \dic_P]$, with $\dic_i \in \R^{M \times P_i}$;
    \item Initialise the weight $\bU^0$ as $\bU^0=[\bU^0_1,\ldots, \bU^0_P]$, $\bU^0_i =[\theta^0_{i1},\ldots, \theta^0_{iP_i}]$, with $\theta^0_{ij}=0$.
                      \end{enumerate}

     \For  {$k=0, \ldots, k_{\max}$}
%            \State Solve the reweighted $\ell_1$-minimisation problem with convex constraints on $\bw$
%            \begin{equation}
%            \min_{\bw}  \frac{1}{2}\| \dic\bw-\by\|_{2}^{2}+\noise\|\bU^{k}\bw\|_1;
%            \end{equation}

    \State Apply Algorithm \ref{alg:splitting} and \ref{alg:weightedlasso} to to get an estimate on $\bw^{k+1}$ to (\ref{weightedlasso});
    %\State Update $\bU^{k+1}$
%     \For  {$k=1, \ldots, k_{\max}$}
%            \State Update $\bw_i^{k+1}=(\dic_i^T\dic_i+\hat{\rho}\bU_i^T\bU_i)^{-1}(\dic_i^T\bb+\hat{\rho}\bU_i^T(\hat{\bz}_i-\hat{\bu}_i))$;
%            \State Update $\hat{\bz}_i^{k+1}=S_{\hat{\lambda}/\hat{\rho}}(\bU_i\bw_i^{k+1}+\hat{\bu}_i^k)$;
%            \State Update
%            $\hat{\bu}_i^{k+1}= \hat{\bu}_i^k+\bU_i\bw_i^{k+1}-\hat{\bz}_i^{k+1}$;
%        \If   {A stopping criterion is satisfied or when $k$ reaches a predefined iteration number $k_{\max}$}
%            \State Break;
%        \EndIf
%     \EndFor
%
%            \State Prune out the estimated zero weights, set new $\dic$;
            \State Set
            $\bU^{(k)}\define \diag\left[\btheta^{(k)}\right]^{-1}$, $\bW^{(k)}\define \diag\left[|\bw^{(k)}|\right];$
%            $$
%            \bU^{k}\define \diag\left[\btheta^{k}\right]^{-1}, \
%           \bW^{k+1}\define \diag\left[|\bw^{k+1}|\right];
%            $$

            \State Update $\theta_j^{(k+1)}$ for the next iteration
            $$
            \theta_j^{(k+1)}=\left[\dic_j^{\mathrm{T}}\left(\noise\mathbf{I}+\dic \bU^{(k)} \bW^{(k+1)} \dic^T\right)^{-1}{\dic_j})\right]^{\frac{1}{2}}; \notag
            $$

        \If   {A stopping criterion is satisfied}
            \State Break;
        \EndIf
     \EndFor
\end{algorithmic}
\end{algorithm}

\section{Numerical Illustration}
\label{sec:example}

\subsection{An Example of Kuramoto Oscillator}

%\begin{example}
%\label{example:k}
A classical example in physics, engineering and biology is the Kuramoto oscillator network (\cite{strogatz2000kuramoto}). We consider a network where the Kuramoto oscillator are nonidentical (each has its own natural oscillation frequency $\omega_i$) and the coupling strengths between nodes are not the same. The corresponding discrete-time dynamics can be described by
\begin{equation}
\begin{aligned}
{\phi_i} (t_{k+1})&={\phi_i} (t_{k})+(t_{k+1}-t_{k})\\
&\left[\omega_i+\sum_{j=1}^{n} w_{ij}g_{ij}(\phi_j(t_k)-\phi_i(t_k))+\xi_i(t_k)\right],
\end{aligned}\label{Kuramoto}
\end{equation}
where $i=1, \ldots, n$, $\phi_i\in[0,2\pi)$ is the phase of oscillator $i$, $\omega_i$ is its natural frequency, and the coupling function $g_{ij}$ is usually taken as sine for all $i, j$. $w_{ij}$ represent the coupling strength between oscillators $i$ and $j$ thus $[w_{ij}]_{n\times n}$ defines the topology of the oscillator network.
Here, assuming we don't know the exact form of $g_{ij}$, we reconstruct from time-series data of the individual phases $\phi_i$ a dynamical network consisting of $n$ Kuramoto oscillators, i.e., we identify the coupling functions $g_{ij}(\cdot)$ as well as the model parameters, i.e., $\omega_i$ and $w_{ij}$, $i, j=1,\dots,n$.

To define the dictionary matrix $\dic$, we assume that all the dictionary functions are functions of a pair of  state variables only and consider $2$ candidate coupling functions $g_{ij}$:  $\sin(x_j-x_i)$, $\cos(x_j-x_i)$. Based on this, we define the dictionary matrix as
\begin{equation}
\begin{aligned}
&\dic_{ij}(x_j(t_k),x_i(t_k)) \define \\
&[\sin(x_j(t_k)-x_i(t_k)), \cos(x_j(t_k)-x_i(t_k))] \in \R^{2}. \notag
\end{aligned}
\end{equation}

To also take into account the natural frequencies, we add to the last column of $\dic_i$ a unit vector. This leads to the following dictionary matrix $\dic_i$:
{\small
\begin{equation}
\begin{aligned}
&\dic_{i} \define \\
&\left[
\begin{array}{cccc}
\dic_{i1}(x_1(t_0),x_i(t_0))& \ldots  &\dic_{in}(x_n(t_0),x_i(t_0)) & 1\\
\vdots  & \vdots  & \vdots  & \vdots \\
\dic_{i1}(x_1(t_{M-1}),x_i(t_{M-1})) & \ldots  & \dic_{in}(x_n(t_{M-1}),x_i(t_{M-1})) & 1
\end{array}
\right] \notag \\
&\in {\mathbb{R}}^{M\times {(2n+1)}}. \notag
\end{aligned}
\end{equation}
}

Then the output can be defined as
$$\by_i\define
\left[\frac{\phi_i(t_{1})-\phi_i(t_{0})}{t_{1}-t_0},\ldots,\frac{\phi_i(t_{M})-\phi_i(t_{M-1})}{t_{M}-t_{M-1}}\right]^{\mathrm{T}}\in {\mathbb{R}}^{M}.$$

To generate the time-series data, we simulated a Kuramoto oscillator network for which $10\%$ of the non-diagonal entries of the weight matrix $[w_{ij}]_{n\times n}$ are nonzero (assuming $g_{ii}$ and $w_{ii}$ are zeros), and the non-zero $w_{ij}$ values are drawn from a standard uniform distribution on the interval $[-10,10]$. The natural frequencies $\omega_{i}$ are drawn from a normal distribution with mean $0$ and variance $10$.
In order to create simulated data, we simulated the discrete-time model \eqref{Kuramoto} and took `measurements data points' every $t_{k+1}-t_{k}=0.1$ between $t = 0$ and $t = 100$ (in arbitrary units) from random initial conditions  which are drawn from a standard uniform distribution on the open interval $(0,2\pi)$.
Thus a total of 1001 measurements for each oscillator phase $\phi_i$ are collected (including the initial value).
Once again, it should be noted that the the number of rows is less than that of columns of the dictionary matrix.

%In the noiseless case, we used Algorithm \ref{alg:weight} with the collected time-series data points per oscillator to reconstruct a network consisting of $100$ randomly coupled Kuramoto oscillators with $10\%$ nonzero links (see Fig.~\ref{fig:K1}), thus $\dic_i\in \R^{450\times 501}$.
%In the reconstruction process, the nonlinearity $\sin(x_j-x_i)$ is the only one selected by our algorithm from the set of candidate functions.
%We further define the weight matrix as
%\begin{eqnarray}
%%{\small
%\bw=
%\left[
%\begin{array}{cccc}
%w_{11} & \ldots & w_{1n} & \omega_1\\
%\vdots & \vdots& \vdots & \vdots \\
%w_{n1} & \ldots & w_{nn} & \omega_n  \notag
%\end{array}
%\right]\in \R^{100\times 101}.
%%}
%\end{eqnarray}
%The true (resp. absolute error) weight matrix, where the last ($101^{th}$) column corresponds to the natural frequencies $\omega_i$ (resp. natural frequency error), is shown in Fig.~\ref{fig:K1} (resp. Fig.~\ref{fig:K2}).

\subsection{Algorithmic Performance Comparisons}

%\guy{``One way is to carry out these updates'': WHICH UPDATES???}

The reweighted lasso algorithm can be implemented in a centralised way by using CVX (\cite{grant2008cvx}) or YALMIP (\cite{lofberg2004yalmip}), Matlab packages for specifying and solving convex optimisation problems. CVX or YALMIP calls generic SDP solvers (SDPT3 (\cite{toh1999sdpt3} or SeDuMi (\cite{sturm1999using}) to solve the problem. While these solvers are reliable for wide classes of optimisation problems, they are not customised for particular problem families, such as ours.

We compare the centralised algorithm using CVX, the centralised algorithm using ADMM,
and the distributed algorithm using ADMM. We fixed the number of measurements $M$ to be $1001$ and varied the network size $n$ between $500$ and $100,000$. For the distributed algorithm, we split the problem into 1000 subproblems where each one has the same dimension. The algorithm is implemented in MATLAB R2012b. The calculations were performed on HP workstation with two 8 core Intel{\circledR} Xeon(R) CPU E5-2650 2.00GHz with 64 GB RAM.

Since the reconstruction problem in (\ref{problem:expand}) for each node is independent, we therefore consider the performance of a single node for illustration. We first investigate the performance for different signal-to-noise ratios of the data generated for this example.
We define signal-to-noise ratio (SNR) by $\mathrm{SNR(dB)}\define 20 \log (\|\dic \bw_{\mathrm{true}}\|_2/\|\bXi\|_2)$. We considered SNR ranging from 5 dB to 25 dB for each generated weight. To characterise the accuracy of a reconstruction, we use the normalised mean square error (NMSE) as a performance index, defined by $\|\hat{\bw}-\bw\|/\|\bw\|$, where $\hat{\bw}$ is the estimate of the true weight $\bw$.
%\guy{Check that this preceding sentence is correct.}
For each SNR, we generated 50 independent experiments (with different initial conditions and parameters) and calculated the average NMSE for each SNR over these 10 experiments. The results are shown in  Fig. \ref{fig:snr}.
\begin{figure*}
\label{fig:snr}
\center
\includegraphics[scale=0.5]{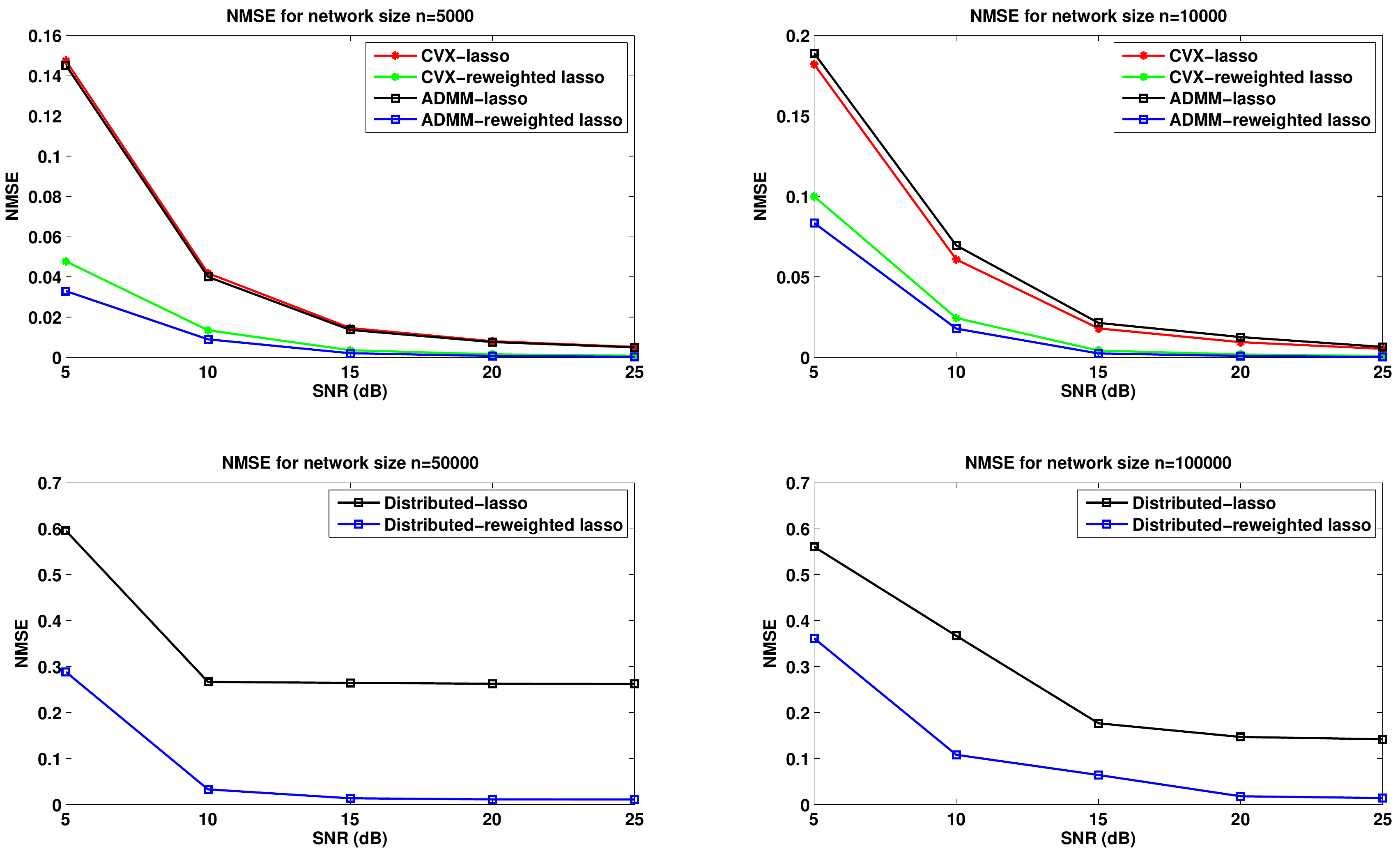}
\caption{Normalised Mean Square Error (NMSE) averaged over 50 independent experiments for the signal-to-noise ratios 5dB, 10 dB, 15 dB, 20 dB, and 25 dB, with different network dimensions $5000$, $10000$, $50000$ and $100000$.}
\label{fig:snr}
\end{figure*}

 %This is no analytical formula to this problem as far as we know.
Next, we compare the computation time for different network sizes for each method. These methods are tested under different SNRs as above. In the implementation of  the distributed algorithm, we use Matlab command \verb"parfor" to parallelise the $\bw_k$-update in Algorithm \ref{alg:splitting}. We use Matlab command \verb"matlabpool(`size')" to start a worker pool. The \verb"size" varies from 2 to 10. For each method, we found that the computation time over each SNR varied slightly (at least within the same magnitude). We calculated the average computation time from a total 250 (=5$\times 50$) independent experiments for each method.
%It should be noted that we calculated the computation time for distributed algorithm by adding the worst computation time for each independent $\bw_i$-update to other centralised computations such as cache of temporary variables, $\bz$-update, $\bu$-update and so on.
The results are shown in Table \ref{time-compare} (unit in second).
For larger problem, with network size 50,000 and 100,000, CVX based reweighted lasso run into memory difficulties. Thus there are no results reported in the table for these network sizes.
For the distributed reweighted lasso algorithm, the computation time decreases when the problem is split between an increasing number of processors.
%Clearly, ADMM based algorithms outperform the generic algorithm. Furthermore, the distributed algorithm reduced. We generated 50 independent experiments (with different initial conditions and parameters) and calculated the average speed for each $N$.
%Indeed, the computation time of ADMM based algorithm can be dramatically reduced if the algorithm is implemented in C instead of Matlab.
However, it should be noted the computation time for distributed reweighted lasso is small partially because Matlab performs some matrix  computations in parallel. On the other hand, CVX exploits only one core.

For all the experiments, we set the penalty parameter $\rho=1$ in the augmented Lagrangian and the scalar regularisation parameter $\lambda=0.05\|\dic^T\by\|_{\infty}$. We also considered termination tolerances $\epsilon^{abs}=10^{-4}$ and $\epsilon^{abs}=10^{-2}$ and set the ADMM iteration number to be 200 and the reweighted iteration number to be 10 throughout all algorithms.

{\large
\begin{table*}
\caption{} \label{time-compare}
%\captionof{table}{Comparison of computation time for different methods and network size}
\begin{center}
\begin{tabular}{|l|*{5}{c|}}\hline
\backslashbox{Methods}{Network Size $n$}
%&\makebox[3em]{5/31}&\makebox[3em]{6/1}&\makebox[3em]{6/2}
%&\makebox[3em]{6/3}&\makebox[3em]{6/4}\\\hline
&$500$&$5000$&$10000$&$50000$&$100000$ \\ \hline
  CVX-reweighted lasso&81.9 & 428.4 & 1218.3& N/A & N/A  \\ \hline
%  ADMM-reweighted lasso&0.6 & 31.8 & 206.6& N/A  & N/A \\ \hline
  Distributed reweighted lasso with 2 cores &24.7 & 84.3 & 156.5   & 587.1   & 1341.5 \\  \hline
  Distributed reweighted lasso with 4 cores &16.4 & 64.1 & 91.5  & 411.8 & 868.7 \\  \hline
  Distributed reweighted lasso with 10 cores &15.5 &  46.9 & 62.9 & 345.2 & 788.3 \\  \hline
\end{tabular}
\end{center}
\center
%\bigskip
%Comparison of computation time for different methods and network size averaged over 50 independent experiments of different SNR.
%\end{center}
\end{table*}
}

\section{CONCLUSION AND DISCUSSION}\label{sec:conclusion}
In this paper, a new distributed reconstruction method for nonlinear dynamical network is proposed. The proposed method only requires time-series data and some prior knowledge about the class of systems for which a dynamical model needs to be built. The network reconstruction problem can be casted as a sparse linear regression problem. Under Bayesian interpretation, this problem is solved using a reweighed $\ell_1$ algorithm which can further reduce the Normalised Mean Square Error in comparison with the classic lasso algorithm. Furthermore, our distributed algorithm can deal with  networks comprising more than 50,000 nodes, which centralised algorithm typically cannot deal with.

Although the convergence of ADMM algorithm and reweighted lasso algorithm have been studied previously (\cite{boyd2011distributed, wipf2007new}), the convergence of the algorithms developed in this paper still need to be properly characterised. We are currently establishing such convergence results as well as the associated convergence rates and their dependence on parameters such as $\rho$ and $\lambda$.

\section{ACKNOWLEDGEMENT}
Authors would like to thank Dr Ye Yuan and Dr Jorge Gon\c{c}alves for helpful discussions and suggestions.
\bibliography{ifacconf}
\end{document}